\documentclass{amsart}

\begingroup
\swapnumbers
\theoremstyle{plain}

\newtheorem{thm}{Theorem}[section]
\newtheorem{cor}[thm]{Corollary}
\newtheorem{prop}[thm]{Proposition}
\newtheorem{lem}[thm]{Lemma}

\theoremstyle{definition}
\newtheorem{dfn}[thm]{Definition}
\newtheorem{rmk}[thm]{Remark}
\newtheorem{rmks}[thm]{Remarks}
\newtheorem{egz}[thm]{Example}

\newtheorem{gen}[thm]{}
\endgroup

\numberwithin{equation}{section}

\DeclareMathOperator{\cDom}{\mathcal{C}\mathcal{D}\mathcal{O}\mathcal{M}}

\newcommand {\Bh}{\mathcal{B}(H)}  
\newcommand {\Bx}{\mathcal{B}(X)}  
\newcommand {\B}{\mathcal{B}}  
\newcommand {\ld}{\lambda} 
\newcommand{\myqed}{\ \hfill \rule{2mm}{2mm}} 

\newcommand {\N}{\mathbb{Z}_{+}} 
\newcommand {\R}{\mathbb{R}} 
\newcommand {\C}{\mathbb{C}} 
\newcommand {\D}{\mathbb{D}} 

\def\firsta{\  {\text {\rm (a)}}\ \ }

\def\condb{\smallskip  \noindent {\text {\rm (b)}}\ \ }

\newcommand {\noi}{\smallskip \noindent} 

\begin{document}

\title[Completely polynomially dominated operators]
{Operators near completely polynomially dominated ones and similarity 
problems}

\author{C. Badea}
\address{D\'epartement de Math\'ematiques, UMR 8524 au CNRS, 
Universit\'e de Lille I, 
F--59655 Villeneuve d'Ascq, France}
\email{Catalin.Badea@agat.univ-lille1.fr}
\urladdr{http://www-gat.univ-lille1.fr/\~{}badea/}

\begin{abstract}  
Let $T$ and $C$ be two Hilbert space operators. 
We prove that if $T$ is near, in a certain sense, 
to an operator completely polynomially dominated with a finite bound 
by $C$, then $T$ is similar to an operator which is completely 
polynomially dominated by the direct sum of $C$ and a suitable 
weighted unilateral shift.  
Among the applications, a refined Banach space version of Rota similarity 
theorem is given and partial answers to a problem of K. Davidson and 
V. Paulsen are obtained. The latter problem 
concerns CAR-valued Foguel-Hankel operators which are generalizations of 
the operator considered by G. Pisier in his example 
of a polynomial bounded operator not similar to a contraction. 
\end{abstract}
\date{}
\maketitle

\bigskip

\section{Introduction}
\thispagestyle{empty}
\begin{gen}{\bf Preamble.}
A good part of the literature concerning similarity problems for 
operators on a Hilbert space was motivated by a single problem. This 
problem asks for a simple criterion to
determine whether a Hilbert space operator is similar 
to a contraction.
The corresponding problems for similarity to 
isometries or unitaries have been solved at the late 1940's by 
Sz. Nagy \cite{sznagy:1947}. The conjectured \cite{halmos:bams} 
characterization~: "an operator 
is similar to a contraction if and only if it is polynomially 
bounded"
was recently shown to be false by G. Pisier \cite{pisier:jams}. Recall 
that $T$ is said to be \emph{polynomially bounded} if there exists a 
constant $M$ such that 
\begin{equation}
	\|p(T)\| \leq M\sup\{ |p(z)| : |z| = 1\}
	\label{eq:01}
\end{equation}
for all polynomials $p$. We refer to \cite{davidson:survey} for the 
history of this counterexample.

A positive answer for the similarity problem was given in 
\cite{paulsen:criterion}.  
The quantitative criterion of V.~Paulsen \cite{paulsen:criterion} 
asserts that an operator 
$T$ is similar to a 
contraction if and only if $T$ is a \emph{completely polynomially bounded} 
operator, which means that equation (\ref{eq:01}) holds 
for all matrix-valued polynomials. 
Moreover, the similarity constant coincides with the 
smallest possible constant $M$ in the analogue of (\ref{eq:01}). 
A more general result for similarity of algebra homomorphisms 
to completely contractive 
ones was proved in \cite {paulsen:pams} (cf. also \cite{paulsen:book}).

Paulsen's criteria are consistent with a variety of similarity 
results 
in operator theory. They are also consistent with results in some areas 
of operator 
algebras and operator spaces theory, areas where completely positive  
and completely bounded maps 
have found to be central tools. Generalizations to Banach space 
operators and to $p$-complete bounded homomorphisms are given in 
\cite{pisier:indiana} (see also \cite{pisier:book}).

We introduce in this paper the notion of operators $T$ (completely) 
polynomially dominated with finite bound by a given operator $C$. 
For instance, we will say that $T$ is \emph{polynomially dominated} with finite
bound by $C$ 
if there exists $M > 0$ such that 
\begin{equation}
	\|p(T)\| \leq M\|p(C)\|
\end{equation}
for all polynomials $p$. 
Completely polynomially dominated operators with finite bound 
generalizes completely polynomially bounded operators. 

The main goal of this note is to show that an operator $T$ near, in a 
certain sense, to a Hilbert space operator completely polynomially 
dominated with a finite bound by a given operator $C$ is similar to 
an operator which is completely polynomially dominated by the direct 
sum of $C$ with a suitable weighted unilateral shift. 
The nearness 
condition for Hilbert space operators (called here 
$\beta$-quadratic nearness) is defined 
in Section 2.
In particular, the class of operators similar to contractions is 
stable under quadratic nearness. A precursor of results of this 
type is \cite{holbrook:sz}.

Applications to similarity problems for Hilbert space 
operators include two partial results concerning an open 
question \cite{davidson/paulsen} of K. Davidson and V. Paulsen. 
The question mentioned in 
\cite{davidson/paulsen} asks for a characterization of those square 
summable sequences for which the corresponding CAR-valued 
Foguel-Hankel operators are similar to contractions. Note that the 
counterexamples of Pisier \cite{pisier:jams} are 
operators of this type. It was this question which was the starting 
point of this note.

Even if the emphasis here will be on Hilbert space operators, 
we will also consider Banach space operators in Theorem 
\ref{thm:main1banach}. As an application, a refined version of Rota's 
\cite{rota} similarity result will be obtained. We will show that, 
given $p > 1$ and a Banach space operator $T$ on $X$ with spectral 
radius less than one, $T$ is similar to an operator $T_1$ on a Banach 
space which, in some sense, "looks like $X$" such that $T_1$ is completely 
polynomially dominated by the unilateral shift $S$ on $\ell_p(X)$. 
This is related to a conjecture of V.I. Matsaev concerning 
contractions on $L_{p}$-spaces.

We also consider the (easiest) corresponding problem for operators near 
ones which are
similar to unitaries or isometries.
We prove that operators asymptotically near operators similar to 
unitaries/isometries are themselvs similar to unitaries or 
isometries. There are polynomially bounded operators which are 
asymptotically near to a contraction without being similar to a 
contraction.

\noi {\bf Acknowledgment.}
Parts of the present paper were written while the author attended the 
Semester on Operator Spaces and Free Probability at Institut 
Henri Poincar\'e, Paris, 1999-2000.
I want to thank L. Kerchy, C. Le Merdy and V. Paulsen for useful 
discussions, suggestions and some simplifications of the arguments 
in an early version. 
\end{gen}

\begin{gen}{\bf Organization of the paper.}
After this preamble we recall some notation, definitions and 
known results. We introduce in the next section the notions of completely
polynomially dominated operators and of
asymptotically near and quadratically near operators. The main 
results in the Hilbert space situation are stated in Section 3. 
This section also 
contains an example of a polynomially bounded operator which is 
asymptotically near to a contraction without being similar to a contraction. 
In Section 4 
the proof of Theorem \ref{thm:main}
is reduced to the proof of Theorem \ref{thm:main1*}. A more 
general version of Corollary  \ref{thm:main1} is stated in the Banach space 
context (Theorem \ref{thm:main1banach}). 
Section 5 
contains several applications to operators similar to contractions, 
including a sufficient condition for the similarity to contractions 
of some CAR-valued Foguel-Hankel operators (Corollary \ref{bthree}) and a 
Banach space Rota theorem (Corollary \ref{thm:rotabanach}).  The proof of 
Theorem \ref{thm:main1banach} is given in Section 6 while the last 
Section contains proofs of the remaining results. 
\end{gen}

\begin{gen}{\bf Preliminaries.}
We recall now some definitions and results and introduce some 
notation. 
We refer to \cite{pisier:book} and \cite{paulsen:book} for more 
information.

\noi {\bf General notation.} By $H, K$ (and $X, Y, E$), 
with or without subcripts, we 
will designate complex Hilbert (respectivelly Banach) spaces. We denote by 
$\Bx$ the algebra of all bounded linear operators on $X$. By operator 
we always 
mean a bounded linear operator. The adjoint of a Hilbert space 
operator $T$ is denoted by $T^{\ast}$. 

\noi {\bf Similarity.} Two Hilbert 
space operators $T_{1}, T_{2} \in \Bh$ are called \emph{similar} if 
there exists an invertible operator $L \in \Bh$ such that $T_{2} = 
L^{-1}T_{1}L$. 

If $\mathcal{A}$ is a class of bounded linear operators, then the 
\emph{similarity constant} $C_{sim}(T_{1},\mathcal{A})$ of $T_1$ with respect to
$\mathcal{A}$ is defined by
$$C_{sim}(T_{1},\mathcal{A}) = \inf \{\|L^{-1}\|\cdot\|L\| : L \in \Bh ,
L^{-1}T_1L \in \mathcal{A} \} .$$

We recall that $T\in \Bh$ is similar to a contraction if and only if 
there exists a \emph{Hilbertian}, equivalent norm on $H$ 
with respect to which $T$ is a contraction.  

\noi {\bf Completely bounded maps.} Let $\mathcal{S} \subset \Bh$ be 
a 
subspace. Let $\varphi : \mathcal{S} \to \mathcal{B}(K)$ be a linear 
map.
Let $M_{n}(\mathcal{S})$ and $M_{n}(\mathcal{B}(K))$ be the spaces of 
matrices with entries respectively in $\mathcal{S}$ and 
$\mathcal{B}(K)$. 
Endow them with the norm induced respectively by $\B(\ell_{n}^2(H))$ 
and $\B(\ell_{n}^2(K))$. The map  $\varphi$ is called 
\emph{completely 
bounded} if there is a constant $M$ such that
$$\sup_{n}\|I_{M_{n}}\otimes \varphi : M_{n}(\mathcal{S}) \to 
M_{n}(\mathcal{B}(K))\| \leq M . $$
The completely bounded (cb) norm $\|\varphi\|_{cb}$ is the smallest 
constant $M$
for which this holds. We call $\varphi$ \emph{completely contractive} if 
$\|\varphi\|_{cb} \leq 1$. The map $\varphi$ is \emph{completely positive} 
if $I_{M_{n}}\otimes \varphi$ is a positive map for each $n$.

The following (Wittstock-Paulsen-Haagerup) 
factorization theorem for completely 
	bounded maps holds \cite[Ch. 3]{pisier:book}, 
\cite[Ch.7]{paulsen:book}~: 
	If $\mathcal{S} \subset \Bh$ is a 
subspace and $\varphi : \mathcal{S} \to \mathcal{B}(K)$ is a linear 
completely bounded map, then there exist a Hilbert space $H_{\pi}$, 
	a unital C$^{\ast}$-algebraic representation $\pi : \B (H) \to 
	\B(H_{\pi})$ and operators $V_2 : K \to H_{\pi}$, $V_1 : H_{\pi}\to 
K$, 
	with $\|V_{1}\| 
	\|V_{2}\| \leq \|\varphi\|_{cb}$, such that $\varphi(a) = 
	V_{1}\pi(a)V_{2}$ for any $a\in \mathcal{S}$.

Let $A(\D)$ be the disk algebra. For an operator $T$, let $\Phi_{T}$ 
be the functional calculus map $p\to p(T)$ defined on polynomials. 
Then $T$ is completely 
polynomially 
bounded if and only if $\Phi_{T}$ extends to a completely bounded map 
on $A(\D)$, if and only 
if $T$ is similar to a contraction \cite{paulsen:criterion}.

Let $p \geq 1$. Similar notions of $p$-complete bounded maps are 
defined in the Banach space context \cite{pisier:book}. If 
$\mathcal{S} \subset \Bx$ is a subspace, a linear map  $\varphi : 
\mathcal{S} \to \mathcal{B}(Y)$ is $p$-completely bounded if 
$$\|\varphi\|_{pcb} := \sup_n \|I_{\B(\ell_p^n)}\otimes \varphi : 
M_n(\mathcal{S}) \to M_n(\B(Y))\| < +\infty ,$$
where $M_n(\B(Y))$ and $M_n(\mathcal{S})$ are now equipped with the 
norms induced by $\B(\ell_p^n(Y))$ and respectively $\B(\ell_p^n(X))$.

We refer to \cite{pisier:indiana}, \cite{pisier:book} for more on 
this, including a factorization theorem.

\noi {\bf Banach spaces of class $SQ_p$.}
Let $p\geq 1$ be a real number. A Banach space $E$ is said to be a 
$SQ_p$-space if it is a quotient of a subspace of an $L_p$-space. 

Let $X$ be a Banach space. A Banach space  $E$ is said to be a 
$SQ_p(X)$-space if it is (isometric to) a quotient of a subspace of 
an ultraproduct of spaces of the form $L_p(\Omega,\mu,X)$. 
Since ultraproducts of $L_p$-spaces is an $L_p$-space, the latter
definition is consistent with the former. The case $p = 2$ 
corresponds to the Hilbertian situation.

$SQ_p(X)$-spaces are characterized by a theorem of Hernandez \cite{hernandez}. 
See also \cite{pisier:indiana} for a different proof using 
$p$-completely bounded maps. Namely, $E$ is a $SQ_p(X)$-space 
if and only if 
$$\|a\|_{p,E} \leq \|a\|_{p,X}$$
for each $n\geq 1$ and each matrix  $a = [a_{ij}]\in M_{n}(\C)$. Here 
$$\|[a_{ij}]\|_{p,Y} =  \sup 
\left[\left(\sum_i\|\sum_ja_{ij}y_j\|^p\right)^{1/p}\right] ,$$
where the supremum runs over all $n$-tuples $(y_1, \cdots , y_n)$ in 
$Y$ which satisfy $\sum\|y_j\|^p \leq 1$.
\end{gen}

\noi {\bf CAR-valued Foguel-Hankel operators.}
A polynomial bounded operator which is not 
completely polynomially bounded was found in 1997 by G. Pisier 
\cite{pisier:jams}. The 
counterexample was a CAR-valued Foguel-Hankel type operator 
(sometimes called a CAR-valued Foias-Williams-Peller type 
operator). 

To be more specific, let $\Lambda$ be a function from an infinite 
dimensional Hilbert space $H$ into
$\Bh$ satisfying 
the \textit{canonical anticommutation relations}~: for all $u,v\in H$,
\begin{gather*}\label{CAR}
 \Lambda(u)\Lambda(v)+\Lambda(v)\Lambda(u)=0  \\ 
\intertext{and}
 \Lambda(u)\Lambda(v)^*+\Lambda(v)^*\Lambda(u)= (u,v) I .
\end{gather*} 
The range  of $\Lambda$ is isometric to Hilbert space. 
Let $\{e_{n}\}_{n\geq 0}$ be an orthonormal basis for $H$, and let 
$C_n=\Lambda(e_n)$ for $n\geq 0$. For an
arbitrary sequence $\alpha=(\alpha_0,\alpha_1,\dots)$ in $\ell^2$, let
$$
 Y_\alpha = \Bigl[ \alpha_{i+j}C_{i+j} \Bigr] 
$$
be a CAR-valued Hankel operator and 
$$
R(Y_\alpha) =
   \begin{bmatrix}S^{*(\infty)}&Y_\alpha\\0&S^{(\infty)}\end{bmatrix} 
.
$$
be the
corresponding Foguel-Hankel operator  \cite{pisier:jams}, 
\cite{davidson/paulsen}.  
Here $S^{(\infty)}$ is the unilateral forward shift of multiplicity 
$\dim H$. The particular choice of
$\alpha$ made by Pisier was $\alpha_{2^k-1}=1$ for $k\ge0$ and
$\alpha_i=0$ otherwise. In this case $R(Y_\alpha)$ is polynomially 
bounded but not completely polynomially bounded. The following more 
general result holds \cite{pisier:jams}, \cite{davidson/paulsen}~:

\begin{thm}[Pisier, Davidson-Paulsen]\label{thm:dav/paul}
Let $\alpha=(\alpha_0,\alpha_1,\dots)$ be a sequence in $\ell^2$ and 
set 
$$A = \sup_{k\ge0} (k+1)^2 \sum_{i\ge k}|\alpha_i|^2$$
and 
$$B_2 = \sum_{k\ge0} (k+1)^2 |\alpha_k|^2 .$$
The operator $R(Y_{\alpha})$ is polynomially bounded if and only if 
$A$ 
is finite. If $R(Y_{\alpha})$ is similar 
to a 
contraction, then $B_2$ is finite.
\end{thm}
It is an open problem if $B_2$ finite implies $R(Y_{\alpha})$ similar 
to a contraction. A partial answer will be proved in Corollary \ref{bthree}.

\section{Dominance and nearness}
\noindent{\bf Dominance.} We start with several definitions.
	
\begin{gen}{\bf Completely polynomially dominated operators.}
Let $T_{1}$ and $T_{2}$ 
be two Hilbert space operators, not necessarily acting on the same 
space. 
We say that 
$T_{1}$ is \emph{completely polynomially dominated} by $T_{2}$ if
	$$\| \left[ p_{ij}(T_{1}) \right]_{1\leq i,j\leq n}\| 
	\leq \|\left[ p_{ij}(T_{2}) \right]_{1\leq i,j\leq n}\|,$$
	for all positive integers $n$ and all $n\times n$ matrices $\left[ 
	p_{ij}\right]_{1\leq i,j\leq n}$ with polynomial entries. Recall 
	that $\left[ p_{ij}(T)\right]_{1\leq i,j\leq n}$ is identified with 
	an operator acting on the direct sum of $n$ copies of the 
	corresponding Hilbert space in a natural way.  
	Let $\cDom(T)$ be the class of all Hilbert space operators 
	completely polynomially dominated 
	by $T$. 
Let $M > 0$ be a positive constant.
We say that 
$T_{1}$ is \emph{completely polynomially dominated with bound} $M$ by 
$T_{2}$ if
	$$\| \left[ p_{ij}(T_{1}) \right]_{1\leq i,j\leq n}\| 
	\leq M\|\left[ p_{ij}(T_{2}) \right]_{1\leq i,j\leq n}\|,$$
	for all positive integers $n$ and all $n\times n$ matrices $\left[ 
	p_{ij}\right]_{1\leq i,j\leq n}$ with polynomial entries. We say 
that $T_1$ is \emph{completely polynomially dominated with finite 
bound}  by $T_{2}$ if it is completely polynomially dominated with 
bound $M$ for a suitable $M$.
 The least bound of complete dominance of $T_1$ by $T_2$ is denoted 
by $M_{cd}(T_{1},T_{2})$. It is the cb norm of the complete bounded 
map $p(T_2) \to p(T_1)$, $p \in \C[z]$.
	
	Similar notions can be defined in the Banach space context. For 
	instance, we say that $T_{1}\in \B (X_{1})$ is $p$-completely 
	dominated with finite bound 
	by $T_{2}\in \B (X_{2})$ if the map  $p(T_2) \to p(T_1)$,  
$p \in \C[z]$, is $p$-completely bounded.
\end{gen}

\begin{egz}
 The following example gives a (generic) class of completely 
dominated 
 operators. Recall the following useful result \cite{sarason:sz}. 
 Let $H$ be a closed 
subspace of $K$ and let $T = P_{H}R\mid H$, $T\in\B (H)$, be the compression of 
$R \in \B (K)$ to $H$. Here $P_{H}$ is the projection onto $H$. 
Then $R$ is a dilation of $T$ (that is, $T^n = 
P_{H}R^n\mid H$ for all $n$) if and only if 
the subspace $H$ is 
semi-invariant for $R$, that is $H = H_{1}\ominus H_{2}$
for two invariant subspaces $H_{1}$ and $H_{2}$ of $R$.

Let $T_{2}\in \B (H_{2})$ be a Hilbert space operator and let $\pi : \B 
(H_{2}) \to 
\B (H_{\pi})$ be a unital  C$^{\ast}$-representation. Let $H_{1}$ be a 
semi-invariant 
subspace for $\pi(T_{2})$. Let $T_{1}\in \B(H_{1})$ be the compression of 
$\pi(T_{2})$ 
on $H_{1}$. Then $T_{1}$ is completely polynomially dominated by $T_{2}$ 
since $\pi$ is completely contractive.
\end{egz}

The following theorem identifies Hilbert space completely 
polynomially 
dominated operators with finite bound.

\begin{thm} \label{thm:cd}
  A Hilbert space operator $T_{1}$ is completely polynomially 
  dominated by $T_{2}$ if and only if $T_{1}$ is unitarily equivalent 
  to the compression of an operator $R_{2}$ to a semi-invariant 
  subspace, $R_{2}$ being the image of $T_{2}$ by a 
unital C$^{\ast}$-representation. 
  A Hilbert space operator $T_{1}$ is completely polynomially 
  dominated by $T_{2}$ with finite bound if and only if $T_{1}$ is 
  similar to an operator completely polynomially dominated by $T_{2}$ 
  and the similarity constant is the least possible bound of 
dominance.
\end{thm}

\noindent {\bf Proof.} Suppose that $T_{1}\in \B (H_{1})$ 
is completely polynomially 
  dominated by $T_{2}$.
Let $\varphi$ be the linear map defined on the 
subspace of the polynomials of $T_{2}\in \B(H_{2})$ by 
$$\varphi(p(T_{2})) = p(T_{1}).$$
The relation of completely polynomially dominance shows that 
$\varphi$  
is well-defined, unital and completely contractive. Then by Arveson's 
theorem \cite{paulsen:book}, Corollary 6.6, $\varphi$ has an 
extension 
$\tilde{\varphi} : \B (H_{2}) \to \B (H_{1})$ which is a unital 
completely 
positive map. By Stinespring's theorem \cite{paulsen:book}, Theorem 
4.1, there are a Hilbert space $K_{1}$, an isometry $V : H_{1} \to 
K_{1}$ and a unital C$^{\ast}$-representation $\pi : \B(H_{1})\to \B(K_{1})$ 
such that 
$$\tilde{\varphi} = V^{\ast}\pi V .$$
Denote $R_{2} = \pi (T_{2})$. We obtain 
$$T_{1}^n = \tilde{\varphi}(T_{2}^n) = V^{\ast}R_{2}^nV$$ 
for each $n \geq 0$ and so \cite{sarason:sz} $T_{1}$ is 
unitarily equivalent to the compression 
of $R_{2}$ to a semi-invariant subspace. 

If $T_{1}$ is completely polynomially 
  dominated by $T_{2}$ with finite bound, then $\varphi$ is 
completely 
  bounded and, by Paulsen similarity theorem, 
  \cite{paulsen:book}, Theorem 8.1, $\varphi$ is similar to a 
completely 
  contractive map with the similarity constant given by the complete 
  bounded norm of $\varphi$. \myqed

\medskip

Using Paulsen's criterion, 
$T\in \Bh$ is completely polynomially bounded (i.~e.~similar to a 
contraction) 
whenever $T$ is completely polynomially dominated with finite 
bound by a given contraction.
	
\noi {\bf Nearness.}
We introduce the following definitions of nearness which will be used 
in the statement of the main results.

\begin{dfn}
	Two operators $T_{1}$ and $T_{2}$ acting on the same space 
	are said to be 
\emph{asymptotically near} 
	if 
	$$\lim_{n\to \infty}\|T_{1}^n - T_{2}^n\| = 0.$$
\end{dfn}

\begin{dfn}
Let $\beta : \N \to \R_{+}^{\ast}$. Two operators $T_{1}$ and $T_{2}$ 
are 
said to be $\beta$-\emph{quadratically near} if 
$$ 
s := \left[\sup_{N\geq 0}\|\sum_{n=0}^{N} \frac{1}{\beta(n)^2} (T_{1}^n - 
		T_{2}^n)(T_{1}^n - 
		T_{2}^n)^{\ast}\|\right]^{1/2}< +\infty .
$$
The two operators are simply called \emph{quadratically near} if 
this condition holds with $\beta(n) = 1$ for each $n$.
		
We denote $s$ in the above definition by 
$near_{2}(T_{1},T_{2},\beta)$. If 
$\beta(n) = 1$ for each $n$, we call $s$ the \emph{nearness} (or 
the $2$-nearness)
between $T_{1}$ and $T_{2}$.
\end{dfn}

The above definition of $\beta$-quadratic nearness uses the row 
Hilbert space operator structure \cite{pisier:bookIHP}. The following 
result gives an equivalent definition.
\begin{lem}\label{lema}
Let $\beta : \N \to \R_{+}^{\ast}$. $T_{1}$ and $T_{2}$ are 
$\beta$-quadratically near with $near_{2}(T_{1},T_{2},\beta) \leq s$ 
if and only if 
\begin{equation}\label{eq:21}
\sum_{n=0}^{+\infty} 
\frac{1}{\beta(n)^2}\| (T_{1}^n - 
		T_{2}^n)^{\ast}y\|^2 \leq s^2 \| y\|^2 \quad  \quad (y \in H).
\end{equation}
If 
\begin{equation}\label{eq:22}
\sum_{n=0}^{+\infty} \frac{1}{\beta(n)^2}\| T_{1}^n - 
		T_{2}^n\|^2 = u^2 < +\infty ,
\end{equation}
then $T_{1}$ and $T_{2}$ are 
$\beta$-quadratically near with $near_{2}(T_{1},T_{2},\beta) \leq u$.
\end{lem}

\noindent {\bf Proof.} For $N \geq 0$ set 
$$
A_{N} = \sum_{n=0}^{N} \frac{1}{\beta(n)^2} (T_{1}^n - T_{2}^n)(T_{1}^n - 
		T_{2}^n)^{\ast} .
$$
Then $T_{1}$ and $T_{2}$ are 
$\beta$-quadratically near with $near_{2}(T_{1},T_{2},\beta) \leq s$ 
if and only if $\sup_{N}\|A_{N}\| \leq s^2$. 
On the other hand, inequality (\ref{eq:21}) holds if and only if 
$\sup_{N}\omega(A_{N}) \leq s^2$, where  
$$\omega(A) = \sup \{ | \langle Ax|x\rangle | : x \in H, \| x \| = 1 
\}$$
is the \emph{numerical radius} of $A$. The stated equivalence follows 
now from the known fact that $\omega(A) = \|A\|$ for normal 
operators $A$.

The second part follows from the fact that (\ref{eq:22}) implies 
(\ref{eq:21}). 
\myqed

\section{Main results : the Hilbert space case}
The classes of operators similar to isometries or unitaries are stable 
under a common nearness condition.

\begin{prop}\label{prop:iso}
A Hilbert space operator asymptotically near an operator similar to 
an 
isometry (or a unitary) is similar to an isometry (respectively a 
unitary).
\end{prop}

The following example, build upon work by Pisier and Davidson and 
Paulsen, shows that there is a polynomially bounded operator which is 
asymptotically near to a contraction without being 
similar to a contraction. 

\begin{egz}\label{egz:32}
We use the notation recalled in Introduction. Let $(\alpha_{k})$ be 
the sequence in $\ell^2$ given by 
$$\alpha_{k} = (k+1)^{-3/2}(\log (k+1))^{-1/2}, \quad k \geq 0 .$$
Then $\sum_{k\ge0} (k+1)^2 
|\alpha_k|^2$ diverges and thus $R(Y_\alpha)$ is not similar to a 
contraction (cf.~Theorem \ref{thm:dav/paul}). 

On the other hand, for $k > 1$, we have 
\begin{eqnarray*}
	\sum_{i\ge k}|\alpha_i|^2  & \leq & \int_{k}^{\infty}\frac{1}{t^3\log 
t}\, dt  \\
	 & \leq & \frac{1}{\log (k)}\int_{k}^{\infty}\frac{1}{t^3}\, 
dt   \\
	 & \leq & \frac{1}{2\log (k)}\frac{1}{(k+1)^2} .
\end{eqnarray*}
Therefore
$$\lim_{k\to \infty} (k+1)^2 \sum_{i\ge k}|\alpha_i|^2 = 0$$
which, using results from \cite{davidson/paulsen}, implies that 
$$\lim_{k\to \infty} \| R(Y_{\alpha})^k - R(0)^k\| = 0.$$
Thus $R(Y_\alpha)$ is asymptotically near the contraction $R(0) = 
S^{*(\infty)}\oplus S^{(\infty)}$, 
without being similar to a contraction.
Note also that $R(Y_\alpha)$ is polynomially bounded since quantity $A$ is 
finite for this $(\alpha_{k})$. \myqed
\end{egz}

The right condition of nearness for the class of operators 
similar to contractions 
follows from the following theorem. 

Let $\beta : \N \to \R_{+}^{\ast}$. We denote by 
$S_{w(\beta)}$ the forward weighted shift on 
$\ell_{2}$, $S_{w}e_{n} = 
w_{n}e_{n+1}$, with weights 
$$w(\beta)_{n} = w_{n} = \frac{\beta(n+1)}{\beta(n)} \quad (n \geq 0).$$
Then $S = S_{w(1)}$ is the unilateral forward shift on 
$\ell_{2}$ obtained for $\beta(n) = 1, n \geq 0$.    

\begin{thm} \label{thm:main}
	Let $T, R \in \mathcal{B}(H)$ and $C \in \mathcal{B}(H_{c})$. 
Suppose that 
	$R$ is completely polynomially dominated with finite bound by $C$. 
	Let $M = M_{cd}(R,C)$ be the least possible bound for this dominance.
	Let $\beta : \N \to \R_{+}^{\ast}$ and suppose that 
	$T$ is $\beta$-quadratically near $R$. Let $s = 
	near_{2}(T,R,\beta)$.
Then $T$ is similar to an operator completelly polynomially 
dominated by $C\oplus S_{w(\beta)}$. Moreover, the similarity 
constant satisfies
	$$C_{sim}(T,\cDom(C\oplus S_{w(\beta)})) \leq M + \beta(0)s.$$
\end{thm}
If $\beta(n) = 1$ for each $n$ we obtain the following consequence.

\begin{cor}
Let $T, R \in \Bh$ and $C \in \B (H_{c})$. Suppose that $T$ is 
quadratically near $R$ and that $R$ is completelly polynomially 
dominated with finite bound by $C$. Then $T$ is similar to the 
compression of $\pi (C\oplus S)$ to a semi-invariant subspace, 
where $\pi$ is a unital C$^{\ast}$-representation defined on $\B 
(H_{c}\oplus \ell_{2})$.
\end{cor}

For similarity to contractions we have 

\begin{cor}\label{cor:bigcor}
Suppose $R\in \Bh$ is similar to a contraction. 
Let $T\in \Bh$ and suppose that there 
exists $C > 0$ such that
$$\sum_{n\geq 0} \|(T^n - R^n)x\|^2 \leq C\|x\|^2$$
for each $x \in H$. Then $T$ is similar to a contraction. 
\end{cor}
Indeed, according to Lemma \ref{lema}, $T^{\ast}$ is quadratically near 
$R^{\ast}$. Note also that $T$ is similar to a contraction if 
and only if $T^{\ast}$ is. 

\begin{rmk}
Operators having their spectrum in the open unit disk are 
quadratically near 0 (the null operator). 
Therefore operators with spectral radius smaller than $1$ are similar 
to contractions (Rota's \cite{rota} theorem). 
The relation of quadratic nearness is an 
equivalence relation. It is easy to see that the equivalence class of 
the null operator is the class of all operators having their spectrum 
in the open unit disk.
\end{rmk}

\section{A reduction of Theorem \ref{thm:main} and a Banach space 
extension}
The main result Theorem \ref{thm:main} is a consequence of the 
following
result. It is a generalization of a result of Holbrook \cite{holbrook:sz}.
\begin{thm} \label{thm:main1*}
	Let $T \in \mathcal{B}(H)$ and suppose that there exist Hilbert 
space $K$, 
	operators $V_{2} : H \to K$, $V_1 : K\to H$, $C_{1} \in 
\mathcal{B}(K)$,
	and a function $\beta : \N \to \R_{+}^{\ast}$ such that 
\begin{equation}\label{eq:41}
\sup_{N\geq 0}\|\sum_{n=0}^{N} \frac{1}{\beta(n)^2} (T^n - 
		V_{1}C_{1}^nV_{2})(T^n - 
		V_{1}C_{1}^nV_{2})^{\ast}\| = s^2 < +\infty .
\end{equation}	
	Then $T$ is similar to an operator completely polynomially 
	dominated by $C_{1}\oplus S_{w(\beta)}\in \mathcal{B}(K\oplus 
\ell_{2})$. 
	Moreover, the similarity constant satisfies
	$$C_{sim}(T,\cDom(C_{1}\oplus S_{w(\beta)})) \leq \|V_{1}\| 
	\|V_{2}\| + \beta(0)s.$$
\end{thm}

\begin{rmks}\label{rmks:4}
 \firsta If $s = 0$ in the above Theorem, then $S_{w}$ can be omitted 
 in the direct sum.
	
 \condb  For an \emph{arbitrary} $T$ and any \emph{finite} $N$, 
 there are operators $V_{1}$, 
 $V_{2}$ and $C_{1}$ like in Theorem \ref {thm:main1*} such that $T^n = 
 V_{1}C_{1}^nV_{2}$ for $ n = 0, 1, \ldots, N$ (cf. 
 \cite[p.910]{halmos:bams}).
\end{rmks}

\begin{gen}{\bf Theorem \ref{thm:main1*} implies Theorem 
\ref{thm:main}.}
Suppose that 
	$R$ is completely polynomially dominated with finite bound by 
	$C\in\B (H_{c})$ and 
	let $M = M_{cd}(R,C)$ be the least possible bound for this dominance.
	Let $\mathcal{S} \subset \B (H_{c})$ be the subspace of all 
operators 
	$p(C)$, $p\in \C[z]$. Consider the map $\Phi : \mathcal{S} 
	\to \Bh$ defined by $\Phi (p(C)) = p(R)$. Since $R$ is completely 
	polynomially dominated with finite bound by $C$, the map $\Phi$ is 
	completely bounded with $\Phi(I) = I$. 
	According to the factorization theorem, there is a Hilbert space $K$, 
	a unital C$^{\ast}$-algebraic representation $\pi : \B (H_{c}) \to 
	B(K)$ and operators $V_{2} : H \to K$, $V_{1} : K \to H$ with 
$\|V_{1}\|\|V_{2}\| \leq M$ such that $\Phi(p(C)) = 
	V_{1}\pi(p(C))V_{2}$ for each polynomial $p$.
	Set $C_{1} = \pi(C)$. We obtain
	$$R^n = \Phi(C^n) = V_{1}\pi(C^n))V_{2} = 
	V_{1}C_{1}^nV_{2}$$
	with $\|V_{1}\| \|V_{2}\| \leq M$. Since $\pi$ is completely 
	contractive, Theorem
	\ref{thm:main1*} implies Theorem \ref{thm:main}. \myqed
\end{gen}
We also obtain the follwing result.
\begin{cor} \label{thm:main1}
	Let $T \in \mathcal{B}(H)$ and suppose that there exist Hilbert 
space $K$, 
	operators $V_{2} : H \to K$, $V_1 : K\to H$, $C_{1} \in 
\mathcal{B}(K)$,
	and a function $\beta : \N \to \R_{+}^{\ast}$ such that 
$$\sum_{n=0}^{+\infty} \frac{1}{\beta(n)^2}\| T^n - 
		V_{1}C_{1}^nV_{2}\|^2 = u^2 < +\infty .$$	
	Then $T$ is similar to an operator completely polynomially 
	dominated by $C_{1}\oplus S_{w(\beta)}\in \mathcal{B}(K\oplus 
\ell_{2})$. 
	Moreover, the similarity constant satisfies
	$$C_{sim}(T,\cDom(C_{1}\oplus S_{w(\beta)})) \leq \|V_{1}\| 
	\|V_{2}\| + \beta(0)u.$$
\end{cor}

In fact the following Banach space version of Corollary \ref{thm:main1} 
holds (for simplicity, we 
will not deal with estimates of the similarity constant here).

We introduce some notation. Consider the space 
$\ell_{p}(\beta ,X)$ of elements $z= (z_{0},z_{1}, 
\ldots)$, $z_{k}\in X$, endowed with the norm
$$\|z\|_{\ell_{p}(\beta ,X)} = \left( 
\sum_{k}\beta(k)^p\|z_{k}\|^p\right)^{1/p}.$$
The shift operator $S$ acts on $\ell_p(\beta,X)$ by 
$$S(z_{0},z_{1}, \ldots ) = (0, z_{0},z_{1}, \ldots ).$$

\begin{thm} \label{thm:main1banach}
	Let $p$ and $q$ be real numbers greater than $1$ such that 
	$\frac{1}{p}+\frac{1}{q}=1$. 
	Let $T \in \mathcal{B}(X)$ and suppose that there exist a 
	$SQ_{p}(X)$-space $Y$, 
	operators $V_{1} : Y \to X$, $V_{2} : X \to Y$, and 
	$C_{1} \in \mathcal{B}(Y)$, 
	and a function $\beta : \N \to \R_{+}^{\ast}$ such that 
$$\sum_{n=0}^{+\infty} \frac{1}{\beta(n)^q}\| T^n - 
		V_{1}C_{1}^nV_{2}\|^q = s^q < +\infty .$$	
	Then there is a Banach space $E$ which is a $SQ_{p}(X)$-space and an 
	isomorphism $L : E \to X$ such that, if $T_{1} = L^{-1}TL \in 
	\mathcal{B}(E)$, then $T_{1}$ is $p$-completely polynomially 
	dominated by $C_{1}\oplus S \in \mathcal{B}(E\oplus 
	\ell_{p}(\beta,X))$.
\end{thm}
\begin{rmk}
As was communicated to the author by V. Paulsen, it is possible to 
prove in a different way Corollary \ref{thm:main1} using Theorem 
\ref{thm:cd}. We have chosen to present a direct proof of 
its Banach space version because of the 
applications of Theorem \ref{thm:main1banach} which are of 
independent interest. A Banach space version of Theorem
\ref{thm:main} can be given using Theorem \ref{thm:main1banach} and 
the factorization theorem for $p$-completed bounded maps of Pisier 
\cite{pisier:indiana}, \cite{pisier:book}. We will not develop this 
idea here.
\end{rmk}

\section{Several applications}
We present now briefly several applications of the main results.

\noi {\bf A Banach space Rota theorem.} It has already been mentioned 
that 
Rota's theorem is a consequence of 
Corollary \ref{cor:bigcor}. The following application of Theorem 
\ref{thm:main1banach} is a refined 
Banach space version of Rota theorem.

\begin{cor} \label{thm:rotabanach}
	Let $X$ be a Banach space and suppose that $T \in \Bx$ has a 
spectral 
	radius smaller than $1$. Then, for every $p > 1$, there exist 
	a Banach space $E$ which is a quotient of $\ell_{p}(X)$ and an 
	isomorphism $L : E \to X$ such that, if $T_{1} = L^{-1}TL \in 
	\mathcal{B}(E)$, then
	\begin{equation}\label{matsaev}
	\|p(T_{1})\|_{\B (E)} \leq \|p(S)\|_{\B (\ell_{p}(X))}
	\end{equation}
	for each analytic polynomial $p$ ; even more generally, 
	$$\|[p_{ij}(T_{1})]\|_{\B (\ell_{p}^n(E))} \leq \|[p_{ij}(S)]\|_{\B 
	(\ell_{p}^n(X))}$$
	for all matrices of polynomials.
\end{cor}

Equation (\ref{matsaev}) shows in particular that $T_{1}$ is a 
contraction. It was conjectured in 1966 by V. I. Matsaev (see 
\cite{peller:steklov}) that 
$$\|p(T_{1})\| \leq \|p(S)\|_{\B (\ell_{p})}$$
holds for all 
contractions $T_{1}$ on an infinite dimensional $L_{p}$-space. 
Several partial results are now known \cite{peller:steklov} but the 
conjecture is still open. The above theorem shows that if the 
spectral 
radius $r(T)$ of $T\in \Bx$ is smaller than one, 
then $T$ is similar to an operator 
on a quotient $E$ of $\ell_{p}(X)$ completely polynomially dominated 
by $S$ on $\ell_{p}(X)$.

If we ask only for a $SQ_{p}(X)$-space $E$ and 
not for a quotient of $\ell_{p}(X)$, 
the proof of Corollary \ref{thm:rotabanach} follows easily from 
Theorem 
\ref{thm:main1banach}. Indeed, if $r(T) < 1$, and 
$\frac{1}{p}+\frac{1}{q} = 1$, then 
$$\sum_{n\geq 0}\|T^n\|^q < +\infty$$
and thus Theorem \ref{thm:main1banach} is applicable with $C_{1} = 
0$. We postpone the proof of Corollary \ref{thm:rotabanach} (with $E$ 
a quotient of $\ell_{p}(X)$) to the last
Section.

\noi {\bf Operators of class} $C_{\rho}$.
Let $\rho > 0$. Operators of class $C_{\rho}$ are
defined as operators having $\rho$-dilations :
$T \in \mathcal{B}(H)$ is in $C_{\rho}$ if there
exists a larger Hilbert space $K \supset H$ and a unitary operator $U$
on $K$ such that
$$T^nh = \rho P_{H}U^{n}h , \quad h \in H \; .$$
Thus contractions are operators of class $C_{1}$. An operator $T$ is 
in $C_{2}$ if and only if $\omega(T) \leq 1$. We refer to \cite{SNF} for more 
information on operators of class $C_{\rho}$. 

A more general class of operators can be constructed as follows \cite{racz}.
Let 
$(\rho_{n})_{n\geq 1}$ be a sequence of positive numbers. We say 
that $T\in \Bh$ is of class 
$C_{\rho_{1},\rho_{2},\ldots}$ 
if there exists a larger Hilbert space $K \supset H$ and a unitary 
operator $U$
on $K$ such that
\begin{equation}\label{eq:52}
	T^nh = \rho_{n} P_{H}U^{n}h , \quad h \in H \; ,
\end{equation}
for all $n\geq 1$. The operator $T$ satisfies (\ref{eq:52}) if and only if 
the spectrum of $T$ 
is in the closed unit disc and 
$$\mbox{ Re }\left[ I + \sum_{n\geq 
1}\frac{2\lambda^n}{\rho_{n}}T^n \right] \geq 0 \quad  \quad 
(|\lambda | < 1) .$$

\begin{cor}[R\'{a}cz]
Let $(\rho_{n})_{n\geq 1}$ be a sequence of positive numbers. 
Suppose that there exist $k \geq 1$ and $M > 0$ such that 
$$\sum_{n=1}^{\infty} (\rho_{nk} - M)^2 < \infty .$$ 
Then every operator of class $C_{\rho_{1},\rho_{2}, \ldots}$ is 
similar to 
a contraction.
\end{cor}

For the proof, denote $S = T^k$. 
Then $S^n = \rho_{nk}V^{\ast}U^{nk}V$ with a suitable 
isometry $V$ and a unitary $U$. It follows that
$$\|S^n - MV^{\ast}U^{nk}V\| \leq \|S^n - 
\rho_{nk}V^{\ast}U^{nk}V\| + | \rho_{nk} - M|.$$
Using Theorem \ref{thm:main1*}, with $C_{1} = U^k$, it follows that 
$S = T^k$ is similar to a contraction 
and 
thus $T$ has the same property (cf. \cite{halmos:bams}).

If $M = \rho_{1} = \rho_{2} = \cdots = \rho$, we obtain the following result 
originally proved by Sz.-Nagy and Foias in 
1967. 

\begin{cor}[Sz-Nagy-Foias]
Every operator of class $C_{\rho}$ is similar to a contraction.
\end{cor}

\noi {\bf Completely bounded maps on $z^dA(\D)$.}
Let $d\geq 1$ be an integer and let $z^dA(\D)$ be the non-unital 
subalgebra of the disc algebra $A(\D)$ consisting of all functions 
$f\in A(\D)$ such that $f(0) = f'(0) = \cdots f^{(d-1)}(0) = 0$.

What happens if the inequality of 
complete dominance with finite bound holds only for polynomials in 
$z^dA(\D)$~? We consider for simplification only Hilbert space 
operators. We refer to \cite[p.80]{pisier:book} and to \cite{mascioni} 
for related results in the Banach space situation.

\begin{cor}
Let $T\in \Bh$ and $C\in \B (H_{c})$ be two Hilbert space operators 
such that 
$$\| \left[ p_{ij}(T) \right]_{1\leq i,j\leq n}\| 
	\leq M\|\left[ p_{ij}(C) \right]_{1\leq i,j\leq n}\|,$$
	for all positive integers $n$ and all $n\times n$ matrices of 
	polynomials $p_{ij}$ in $z^dA(\D)$. Then $T$ is similar to an 
	operator completely polynomially dominated by $C\oplus S\in \B 
	(H_{c}\oplus \ell_{2})$. 
\end{cor}

For the proof, note that the map $P(C)\to P(T)$ defined on 
the subspace 
$$\{P(C) : P \in z^dA(\D), P \mbox{ polynomial }\}$$ 
is completely bounded. By the factorization theorem \cite[Theorem 
3.6]{pisier:book}, we can write 
$$P(T) = V_{1}\pi(P(C))V_{2} \quad ; \quad P \in z^dA(\D)$$ 
with suitable operators $V_{1}$, $V_{2}$ 
and a unital C$^{\ast}$-algebraic representation $\pi$ on $\B (H_{c})$. 
Let $C_{1} = \pi(C)$. We obtain
$$T^k = V_{1}C_{1}^kV_{2} \quad ; \quad k \geq d .$$
This shows that $T$ is quadratically near $C_{1}$. The conclusion 
follows now from Corollary \ref{thm:main1}.

\begin{cor}[Paulsen criterion for $z^dA(\D)$]
Let $d\geq 1$. Let $T\in \Bh$ and suppose that
$$\| \left[ p_{ij}(T) \right]_{1\leq i,j\leq n}\| 
	\leq M\sup_{|z|=1}\|\left[ p_{ij}(z) \right]_{1\leq i,j\leq n}\|,$$
	for all positive integers $n$ and all $n\times n$ matrices of 
	polynomials $p_{ij}$ in $z^dA(\D)$. Then $T$ is similar to a 
	contraction.
\end{cor}

\noi {\bf CAR-valued Foguel-Hankel operators.} 
We use notation as above.
\begin{cor}\label{bthree}
Let $\alpha=(\alpha_0,\alpha_1,\dots)$ be a sequence in $\ell^2$ such 
that 
$$B_{3} := \sum_{k\ge0} (k+1)^3 |\alpha_k|^2 < +\infty .$$
Then $R(Y_{\alpha})$ is similar to a contraction.
\end{cor}

\noi{\bf Proof.} Set $R(0) = S^{*(\infty)}\oplus S^{(\infty)}$. 
Using the notations of \cite{davidson/paulsen}, we have 
$$\| R(Y_{\alpha})^n - R(0)^n\| \leq \|\mathcal{Y}_{\alpha}(z^n)\|.$$
It was proved in \cite{davidson/paulsen} that 
$$\|\mathcal{Y}_{\alpha}(z^n)\| \leq (n+1) \left[\sum_{i\ge 
n}|\alpha_i|^2\right]^{1/2} .$$
We obtain
$$\sum_{n\geq 0}\| R(Y_{\alpha})^n - R(0)^n\|^2 \leq \sum_{n\geq 0} 
(n+1)^2 \left[\sum_{i\ge 
n}|\alpha_i|^2\right].$$
By a Abel summation method, the series $\sum_{n\geq 0} 
(n+1)^2 \left[\sum_{i\ge n}|\alpha_i|^2\right]$ is convergent if 
$$\sum_{n\geq 0} \left[ \sum_{0\leq i\leq n} (i+1)^2 \right] 
|\alpha_n|^2$$
it is. It is indeed convergent because of 
our assumption on $B_{3}$. Therefore 
$R(Y_{\alpha})$ is quadratically near the contraction $R(0)$ and 
thus similar to a contraction.
\myqed

We still don't know if 
$B_{2}$ finite implies $R(Y_{\alpha})$ similar to a contraction. 
Nevertheless, the following similarity 
result holds.

\begin{cor}
Let $\alpha=(\alpha_0,\alpha_1,\dots)$ be a sequence in $\ell^2$ such 
that 
$$B_{2} := \sum_{k\ge0} (k+1)^2 |\alpha_k|^2 < +\infty .$$
Then $R(Y_{\alpha})$ is similar to an operator completely 
polynomially 
dominated by $R(0)\oplus D$, where $D\in \B (\ell_{2})$ is the Dirichlet shift, 
i.e. the weighted unilateral 
shift with weights $w_{n} = \sqrt{(n+2)/(n+1)}$.
\end{cor}

Note that $R(0)$ is a contraction while the Dirichlet shift is
expansive ; it is however a $2$-isometry \cite{agler/stankus}, that 
is $I - 2D^{\ast}D + D^{\ast2}D^2 = 0$. 

The proof is similar to the proof of the precedent corollary : if 
$\beta(n) = \sqrt{n+1}$, then 
$$\frac{1}{\beta(n)}\| R(Y_{\alpha})^n - R(0)^n\| \leq 
\sqrt{n+1} \left[\sum_{i\ge 
n}|\alpha_i|^2\right]^{1/2} .$$
This shows that
$$\sum_{n\geq 0}\frac{1}{n+1}\| R(Y_{\alpha})^n - R(0)^n\|^2 \leq 
\sum_{n\geq 0} 
(n+1) \left[\sum_{i\ge 
n}|\alpha_i|^2\right]$$
and the right hand side is convergent if $B_{2}<+\infty$. Apply 
Corollary \ref{thm:main1} with $\beta(n) = \sqrt{n+1}$ and $C_{1} = 
R(0)$.

\begin{rmk}
Corollary \ref{bthree} was obtained as a particular case of a general 
theorem. Using other methods, 
Vern Paulsen and the author improved Corollary \ref{bthree} as 
follows : $R(Y_{\alpha})$ is similar to a contraction 
if there exists 
$\varepsilon > 0$ such that
$$B_{2+\varepsilon} := \sum_{k\ge 0} (k+1)^{2+\varepsilon}|\alpha_k|^2 
< 
+\infty .$$
Details will be given elsewhere \cite{bp}. 
A different sufficient condition for the similarity to 
contractions of 
operator-valued 
Foguel-Hankel operators was given by G. Blower \cite{blower}.
\end{rmk}

\section{Proof of Theorem \ref{thm:main1banach}}

Put, for simplicity, $C_{1} = C$.
Let $\gamma$ be a positive constant. We will chose this constant in 
the proof of Theorem \ref{thm:main1*} in the next section when 
estimating the similarity constant. 

Set
\begin{equation} \label{eq:norm}
\left|x\right|^p = \inf \left\{\gamma^p
\left\|\sum_{n\geq 0}C^nV_{2}x_{n}\right\|^p_{Y} + 
\sum_{n\geq 0}\beta(n)^p\left\|x_{n}\right\|^p : x = \sum_{k\geq 0} T^kx_{k}\right\},
\end{equation}
the $\inf$ being taken over all (finite) decompositions of $x$ as 
sums of powers of $T$ applied to elements of $X$. 

\begin{gen}$|\cdot|$ {\bf is a seminorm.}
Take two decompositions
$$x = \sum_{k = 0}^d T^kx_{k}$$
and
$$y = \sum_{k = 0}^d T^ky_{k}$$
for fixed $x$ and $y$ in $X$. By adding eventually $x_{k} = 0$ or 
$y_{k} = 0$, we may assume that decompositions have the 
same length $d+1$. This will be always used in the sequel without any 
further comment.

Using the triangle inequality $\|a+b\| \leq \|a\|+\|b\|$ 
in $\ell^{d+1}_{p}(X)$ for 
$$a = (\gamma\sum_{n=0}^{d}C^nV_{2}x_{n},\beta(0) x_{0},\beta(1) 
x_{1}, \ldots, 
\beta(p) x_{p})$$
and
$$b = (\gamma\sum_{n=0}^{d}C^nV_{2}y_{n},\beta(0) y_{0},\beta(1) 
y_{1}, \ldots, 
\beta(p) y_{p})$$
and taking the infimum over all representations of $x$ and $y$, we get
$$|x+y| \leq |x| + |y|.$$

The proofs of the inequality $|\ld x| \leq |\ld| |x|$ and its 
converse are left to the reader. 
\end{gen}

\begin{gen}$|\cdot|$ {\bf is an equivalent norm.}
The representation $x = x_{0} + Tx_{1}$ with $x_{0} = x$ and $x_{1} = 
0$, gives 
$$|x|^p \leq \gamma^p\|V_{2}x\|^p + \beta(0)^p\|x\|^p \leq 
(\gamma^p\|V_{2}\|^p + \beta(0)^p)\|x\|^p$$
and therefore 
\begin{equation} \label{eq:12}
	|x| \leq \left[\gamma^p\|V_{2}\|^p + \beta(0)^p\right]^{1/p} \; 
	\|x\|.
\end{equation}

For the converse inequality, suppose that 
$$x = x_{0} + Tx_{1} + \cdots + T^dx_{d}.$$
We have 
\begin{eqnarray*}
\|x \| & = & 
\|\sum_{k=0}^{d}V_{1}C^kV_{2}x_{k} + 
\sum_{k=0}^{d}(T^k - V_{1}C^kV_{2})x_{k}\| \\
 & \leq & \frac{1}{\gamma}\|V_{1}\|  
 \gamma\|\sum_{k=0}^{d} C^kV_{2}x_{k}\| + \sum_{k=0}^{d} 
 \frac{1}{\beta(k)}\|T^k - 
  V_{1}C^kV_{2}\| \beta(k) \|x_{k}\| .
  \end{eqnarray*}
By using the H\"{o}lder inequality, the last 
quantity is less or equal than
$$\left[\frac{1}{\gamma^q}\|V_{1}\|^q +
\sum_{k=0}^{d}\frac{1}{\beta(k)^q}\|T^k-V_{1}
C^kV_{2}\|^q\right]^{1/q}\left[\gamma^p\|\sum_{k=0}^{d}C^k
V_{2}x_{k}\|^p + \sum_{k=0}^{d}\beta(k)^p\|x_{k}\|^p\right]^{1/p}.$$

Taking the infimum over all representations of $x$, we obtain
\begin{equation} \label{eq:13}
\|x\| \leq \left[\frac{\|V_{1}\|^q}{\gamma^q} + s^q\right]^{1/q} |x|.
\end{equation}
Thus $|\cdot|$ is a norm equivalent to the original one and, using 
(\ref {eq:12}) and (\ref{eq:13}), we have
\begin{equation}\label{eq:14}
\left[\frac{\|V_{1}\|^q}{\gamma^q} + s^q\right]^{-1/q}\|x\| 
\leq |x| \leq \left[\gamma^p\|V_{2}\|^p + \beta(0)^p\right]^{1/p} \; 
\|x\|.
\end{equation}
We denote by $E$ the Banach space $X$ endowed with the new norm 
$|\cdot|$.
\end{gen}

\begin{gen}{\bf The Banach space} $E$ {\bf is a} $SQ_{p}(X)${\bf 
-space.} Let $x_{j} \in X$, $j = 1, \cdots , 
n$, with their decompositions 
$$x_{j} = \sum_{k\geq 0}T^kx_{j}^{(k)}.$$
Let $a = [a_{ij}]\in M_{n}(\C)$ be a matrix such that $\|a\|_{p,X} 
\leq 1$.
This means that
\begin{equation}\label{hernandez}
\sum_{i}\|\sum_{j}a_{ij}y_{j}\|^p \leq \sum_{j}\|y_{j}\|^p
\end{equation}
for all $y_{j} \in X$, $j = 1, \cdots , n$.
We will then have
$$\sum_{j=1}^n a_{ij}x_{j} 
= \sum_{k}T^k\left( \sum_{j}a_{ij}x_{j}^{(k)}\right).$$
By Hernandez theorem we have to prove that 
$\|a\|_{p,E} \leq \|a\|_{p,X}$. Recall that $Y$ is a $SQ_{p}(X)$-space. 
We have
\begin{eqnarray*}
	\sum_{i}|\sum_{j}a_{ij}x_{j}|^p &   \\
	\leq & \sum_{i}\left(\gamma^p
	\|\sum_{k}C^kV_{2}(\sum_{j}a_{ij}x_{j}^{(k)})\|^p_{Y} +
	 \sum_{k}\beta(k)^p\|\sum_{j}a_{ij}x_{j}^{(k)}\|^p\right)    \\
	= & 
	\gamma^p\sum_{i}\|\sum_{j}a_{ij}(\sum_{k}C^kV_{2}x_{j}^{(k)})\|^p _{Y} + 
	\sum_{k}\beta(k)^p\sum_{i}\|\sum_{j}a_{ij}x_{j}^{(k)}\|^p   \\
	\leq & \gamma^p\sum_{j}\|\sum_{k}C^kV_{2}x_{j}^{(k)}\|^p_{Y} + 
	\sum_{k}\beta(k)^p\sum_{j}\|x_{j}^{(k)}\|^p    \\
	 & (\mbox{ by using Eq. (\ref{hernandez}) for $X$ and $Y$ })   \\
	= & \sum_{j}\left( \gamma^p \|\sum_{k}C^kV_{2}x_{j}^{(k)}\|^p + 
	\sum_{k}\beta(k)^p \|x_{j}^{(k)}\|^p\right)  . 
\end{eqnarray*}
By taking infimum over all possible decompositions we get
$$\sum_{i}|\sum_{j}a_{ij}x_{j}|^p \leq \sum_{j}|x_{j}|^p$$
and therefore $E = (X,|\cdot|)$ is a $SQ_{p}(X)$-space.
\end{gen}

\begin{gen}{\bf The operator} $T$ {\bf with respect to} $|\cdot|$. 
Let $x$ be decomposed as $x = \sum_{k\geq 0} T^kx_{k}$ and let 
$$P(z) = \sum_{s=0}^da_{s}z^s$$
be a fixed polynomial. 
Then 
$$P(T)x = \sum_{k}T^k\left( \sum_{i+j=k}a_{i}x_{j}\right) $$
is a decomposition of $P(T)x$. We obtain
$$|P(T)x|^p \leq \Sigma_{1} + \Sigma_{2},$$
where the two sums are given by 
$$\Sigma_{1} = \gamma^p\|\sum_{k}C^kV_{2}\left( 
\sum_{i+j=k}a_{i}x_{j}\right)\|^p$$
and 
$$\Sigma_{2} = \sum_{k}\beta(k)^p\|\sum_{i+j=k}a_{i}x_{j}\|^p.$$

\noindent {\bf The first sum.} Since 
$$\sum_{k}C^kV_{2}\left( 
\sum_{i+j=k}a_{i}x_{j}\right) = \sum_{m}a_{m}C^m\left( 
\sum_{n}C^nV_{2}x_{n}\right),$$
we have
$$\Sigma_{1} = \gamma^p\|P(C)\left( \sum_{n}C^nV_{2}x_{n}\right)\|^p 
\leq 
\gamma^p\|P(C)\|_{\B(Y)}^p\|\sum_{n}C^nV_{2}x_{n}\|^p.$$

\noindent {\bf The second sum.}
The shift operator on $\ell_{p}(\beta ,X)$, also denoted by $S$ acts 
by 
$$S(z_{0},z_{1}, \ldots ) = (0, z_{0},z_{1}, \ldots ).$$
Denote $\tilde{x} = (x_{0},x_{1}, \ldots) \in \ell_{p}(\beta ,X)$, 
where $x_{k}$ are the elements occuring in the (finite) decomposition of $x$. 
The $n$th component of 
$P(S)\tilde{x}\in \ell_{p}(\beta ,X)$ is $\sum_{i+j=n}a_{i}x_{j}$ ; 
hence
\begin{eqnarray*}
	\Sigma_{2} & = & \sum_{k}\beta(k)^p\|\sum_{i+j=k}a_{i}x_{j}\|^p  \\
	 & = & \|P(S)\tilde{x}\|^p_{\ell_{p}(\beta ,X)}  \\
	 & \leq & \|P(S)\|^p_{\B(\ell_{p}(\beta ,X))}\left( \sum_{n\geq 
0}\beta(n)^p\|x_{n}\|^p\right).
\end{eqnarray*}

Combining now the estimates for the two sums, we obtain
$$|P(T)x|^p \leq \max (\|P(C)\|^p , 
\|P(S)\|^p_{\B(\ell_{p}(\beta ,X))})\left( 
\gamma^p\|\sum_{n\geq 0}C^nV_{2}x_{n}\|^p + 
\sum_{n\geq 0}\beta(n)^p\|x_{n}\|^p\right).$$
Taking the infimum over all representations of $x$ we get
$$|P(T)x| \leq \max (\|P(C)\|_{\B (Y)} , 
\|P(S)\|_{\B(\ell_{p}(\beta ,X))})|x|.$$
Therefore
$$\|P(T)\|_{B(E)} \leq \max (\|P(C)\|_{\B (Y)} , 
\|P(S)\|_{\B(\ell_{p}(\beta,X))}).$$

In an analogous way it can be proved that
$$\|[P_{ij}(T)]\|_{B(\ell^n_{p}(E))} \leq 
\max (\|[P_{ij}(C)]\|_{\B(\ell^n_{p}(Y))} , 
\|[P_{ij}(S)]\|_{\B(\ell^n_{p}(\beta,X))})$$
for all polynomials with matrix coefficients. We omit the details.
\end{gen}

\section{Remaining proofs}

\noi {\bf Proof of Theorem \ref{thm:main1*}}
Set again $C_{1} = C$. 
Consider the equivalent norm $|\cdot|$ as defined in the previous 
proof ($p = q = 2$, $X = H$ and $\gamma$ to be precised later on).
Since the class of Hilbert spaces is stable 
by taking subspaces, quotients and ultraproducts of spaces of the form 
$L_{2}(\mu;H)$, $E$ is Hilbertian, that is, the 
new norm $|\cdot|$ comes from an inner product. Also, the unilateral 
shift $S$ on $\ell_{2}(\beta)$ is unitarily equivalent to the weighted 
shift $S_{w(\beta)}$ on $\ell_{2}$ \cite{shields}. The other parts of 
the preceding proofs, excepting the inequality corresponding to 
(\ref{eq:13}), are the 
same. The proof of the inequality 
$$\|x\| \leq \left[\frac{\|V_{1}\|^2}{\gamma^2} + s^2\right]^{1/2} 
|x|$$
runs as follows. 

Suppose $x = x_{0} + Tx_{1} + \cdots + T^dx_{d}$.
We have 
\begin{eqnarray*}
\|x \| & = & 
\|\sum_{k=0}^{d}V_{1}C^kV_{2}x_{k} + 
\sum_{k=0}^{d}(T^k - V_{1}C^kV_{2})x_{k}\| \\
 & \leq & \frac{1}{\gamma}\|V_{1}\|  
 \|\sum_{k=0}^{d}\gamma C^kV_{2}x_{k}\| + \|\sum_{k=0}^{d} 
  (T^k - 
  V_{1}C^kV_{2})x_{k}\| .
  \end{eqnarray*}
Let $y\in H$. It follows from Lemma \ref{lema} that 
$$\sum_{n=0}^{+\infty} \frac{1}{\beta(n)^2}\| (T^n - 
		V_{1}C_{1}^nV_{2})^{\ast}y\|^2 \leq s^2\|y \|^2.$$
We obtain 
\begin{eqnarray*}
	| \langle \sum_{k=0}^{d}(T^k - V_{1}C^kV_{2})x_{k} , y\rangle | & = & 
	|\sum_{k=0}^{d}\langle \beta(k)x_{k} , 
	\frac{1}{\beta(k)}(T^k - V_{1}C^kV_{2})^{\ast}y\rangle |  \\
	& \leq & \left[\sum_{k}\beta(k)^2\|x_{k}\|^2)^{1/2}\right] \left[ 
	\sum_{n=0}^{d} \frac{1}{\beta(n)^2}\| (T^n - 
		V_{1}C_{1}^nV_{2})^{\ast}y\|^2\right]^{1/2}   \\
	& \leq & \left[\sum_{k}\beta(k)^2\|x_{k}\|^2)^{1/2}\right]s\|y\| .
\end{eqnarray*}
Therefore 
$$\|\sum_{k=0}^{d} (T^k - V_{1}C^kV_{2})x_{k}\| \leq s\left[ 
\sum_{k}\beta(k)^2\|x_{k}\|^2\right]^{1/2} .$$
Another application of the Cauchy-Schwarz inequality yields
\begin{eqnarray*}
	\|x\|  & \leq & \frac{1}{\gamma}\|V_{1}\|  
 \|\sum_{k=0}^{d}\gamma C^kV_{2}x_{k}\| + s \left[ 
\sum_{k}\beta(k)^2\|x_{k}\|^2\right]^{1/2}  \\
	 & \leq & 
	 \left[\frac{1}{\gamma^2}\|V_{1}\|^2 +s^2\right]^{1/2}\left[ \|\sum_{k}\gamma 
	 C^kV_{2}x_{k}\|^2 + \sum_{k}\beta(k)^2\|x_{k}\|^2\right]^{1/2} .
\end{eqnarray*}
Taking the infimum over all representations of $x$, we obtain
$$\|x\| \leq \left[\frac{\|V_{1}\|^2}{\gamma^2} + s^2\right]^{1/2} 
|x|.$$
This gives the 
similarity statement. 

We prove now the estimate for the similarity constant. 
From Equation (\ref{eq:14})  and the proof given above we have
$$
C_{sim}(T,\cDom(C\oplus S_{w(\beta)})) \leq 
\left[\frac{\|V_{1}\|^2}{\gamma^2} + s^2\right]^{1/2} 
\left[\gamma^2\|V_{2}\|^2 + \beta(0)^2\right]^{1/2}.
$$
By assuming $C = 0$ if necessary, we may assume 
that $V_{2}$ is not the null operator. If $s \neq 0$, 
choose
$$\gamma = 
\left[\frac{\beta(0)\|V_{1}^{\ast}\|}{s\|V_{2}\|}\right]^{1/2}.$$
We then have
$$ C_{sim}(T,\cDom(C\oplus S_{w(\beta)}))^2 \leq (\|V_{1}^{\ast}\| 
	\|V_{2}\| + \beta(0)s)^2.$$ 
	
	If $s = 0$, then $T^n = V_{1}C^nV_{2}$ and thus $T$ is completely 
polynomially dominated by $C$ with bound 
$\|V_{1}\|\cdot\|V_{2}\|$. 
Apply now Theorem \ref{thm:cd}. Note that in this case $S_{w(\beta)}$ 
is absent from the direct sum.
The proof of Theorem \ref{thm:main1*} is now complete. \myqed

\noi {\bf Proof of Corollary \ref{thm:rotabanach}.} The proof of this 
version of Rota theorem is similar to the proof of Theorem 
\ref{thm:main1banach}. Indeed, if $C=0$, then the new norm $|\cdot|$ 
is given by 
$$|x|^p = \inf \{\sum_{n\geq 0}\beta(n)^p\|x_{n}\|^p : x = \sum_{k\geq 
0} T^kx_{k}\},
$$
the $\inf$ being taken over all (finite) decompositions of $x$ as 
sums of powers of $T$ applied to elements of $X$. This is the 
quotient norm of $\ell_{p}(X)/\mbox{ Ker }(\psi)$, where the onto map $\psi$ is 
given by  
$$\ell_{p}(X) \ni (x_{0},x_{1}, \ldots ) \mapsto \psi (x_{0},x_{1}, 
\ldots ) = \sum_{k}T^kx_{k} \in X .$$
Take $E$ to be $X$ with this new norm.
The rest of the proof is the same. \myqed

\noi {\bf Proof of Proposition \ref{prop:iso}.} 
For the first part of the theorem, 
it is sufficient to prove that an operator asymptotically near an 
isometry is similar to an isometry. Indeed, if we suppose that 
$$\lim_{n\to \infty}\|T^n - L^{-1}V^nL\| = 0,$$ 
with $V$ an isometry, then 
$$\|(LTL^{-1})^n - V^n\| = \|L(T^n - L^{-1}V^nL)L^{-1}\|$$
tends to $0$ as $n$ goes to infinity and so we will obtain the 
similarity of $LTL^{-1}$, so of $T$, to an isometry. 

Now, if $T$ is asymptotically near an isometry $V$, then for each
$r\in ]0,1[$ there exists $k\in \N$ such that 
$$\sup_{n\geq k}\|T^n - V^n\| \leq r.$$
Set $R = T^k$ and $W = V^k$ ($W$ is an isometry). We obtain
$$\sup_{m\geq 1}\|R^{m} - W^{m}\| \leq r < 1.$$
This implies that, for each $x$ and each $m \geq 1$,
$$(1-r)\|x\| = \|W^mx\| - r\|x\| \leq \|W^mx\| - \|R^mx - W^mx\|\leq 
\|R^mx\| \leq (1+r)\|x\|.$$
By a theorem of Sz.-Nagy \cite{sznagy:1947}, $R = T^k$ is similar to 
an isometry and this implies \cite[Corollary 4.2]{popescu:michigan} 
that $T$ is 
similar to an isometry.

Suppose now that $T$ is asymptotically near a unitary $U$. By the 
first part of the proof, $T$ is similar to an isometry. Therefore we 
can write $V^{\ast} = 
L^{-1}T^{\ast}L$, with $V$ an isometry, for a suitable invertible 
operator $L$. 
But $T^{\ast}$ is asymptotically near the isometry $U^{\ast}$ and so 
$T^{\ast}$ is similar to an isometry. This implies that $T^{\ast}$ 
and 
$V^{\ast}$ are injective and so the isometry $V$ is also onto. 
Therefore $V$ is unitary and so $T$ is similar to a unitary. \myqed

\end{document}